\documentclass[]{interact}

\usepackage{color,soul}

\usepackage{algorithm}
\usepackage[]{algpseudocode}

\makeatletter
\def\algbackskip{\hskip-\ALG@thistlm}
\makeatother

\usepackage[natbibapa,nodoi]{apacite}
\usepackage{epstopdf}
\usepackage[caption=false]{subfig}

\theoremstyle{plain}
\newtheorem{theorem}{Theorem}[section]

\theoremstyle{definition}
\newtheorem{definition}[theorem]{Definition}

\theoremstyle{remark}

\usepackage{tikz}
\usetikzlibrary{arrows.meta,
                chains,
                positioning,
                shapes.geometric
                }

\newcommand{\R}{\mathbb{R}}

\DeclareMathOperator{\Ima}{Im}
\DeclareMathOperator{\rank}{rank}

\usepackage[hidelinks,breaklinks]{hyperref}       

\begin{document}


\title{Topological Machine Learning for Multivariate Time Series}


\author{
\name{Chengyuan Wu\textsuperscript{1,2}\thanks{CONTACT Chengyuan Wu Email: wuchengyuan@u.nus.edu} and Carol Anne Hargreaves\textsuperscript{1}}
\affil{1.\ Data Analytics Consulting Centre, Department of Statistics and Applied Probability, Faculty of Science, National University of Singapore, Singapore}
\affil{2.\ Institute of High Performance Computing, A*STAR, Singapore}
}


\maketitle

\begin{abstract}
We develop a {method} for analyzing multivariate time series using topological data analysis (TDA) methods. The proposed methodology involves converting the multivariate time series to point cloud data, calculating Wasserstein distances between the persistence diagrams and using the $k$-nearest neighbors algorithm ($k$-NN) for supervised machine learning. Two methods (symmetry-breaking and anchor points) are also introduced to enable TDA to better analyze data with heterogeneous features that are sensitive to translation, rotation, or choice of coordinates. We apply our methods to room occupancy detection based on 5 time-dependent variables (temperature, humidity, light, CO\textsubscript{2} and humidity ratio). Experimental results show that topological methods are effective in predicting room occupancy during a time window. {We also apply our methods to an Activity Recognition dataset and obtained good results.}
\end{abstract}

\begin{keywords}
Topological data analysis; machine learning; artificial intelligence; multivariate time series; room occupancy
\end{keywords}

\section{Introduction}
Topological Data Analysis (TDA) is a relatively new branch of data analysis using techniques from topology to study data \citep{Edelsbrunner2002,zomorodian2012topological,Zomorodian2005}. It has been applied with great success in several fields such as biomolecular chemistry \citep{xia2018multiscale,xia2015multiresolution}, drug design \citep{cang2018integration}, and network analysis \citep{Carstens2013}. Notably, the winners of the Drug Design Data Resource (D3R) Grand Challenge have utilized TDA in their algorithms \citep{nguyen2019mathematical,nguyen2019mathdl}. Topological data analysis can be combined with methods in machine learning (including deep learning) \citep{hofer2017deep,nguyen2019mathematical} as well as statistical methods \citep{bubenik2015statistical}. 

Though the term ``topology'' can be used to refer to a wide array of subjects, the topological tools used in TDA generally refer to algebraic topology \citep{letscher2012persistent,wu2020weighted}, or to be specific persistent homology \citep{Ghrist2008,edelsbrunner2012persistent}. Broadly speaking, persistent homology analyzes the ``shape'' of the data to deduce intrinsic properties of the data. Other prominent tools in TDA include Mapper \citep{singh2007topological,ray2017survey} and discrete Morse theory \citep{forman1998morse,forman2002user,wu2020discrete}. Due to the fact that TDA works quite differently from most other data analysis techniques, it can sometimes detect features that are missed by traditional methods of analysis \citep{nicolau2011topology}.

Traditionally, the strengths of TDA include the fact that it analyzes data in a coordinate-free way \citep{lum2013extracting,offroy2016topological} (independent of the coordinate system chosen), as well as being translation-invariant and rotation-invariant \citep{khasawneh2016chatter,bonis2016persistence}. As a direct consequence of these strengths, however, it may be hard for TDA to effectively analyze data that is sensitive to choice of coordinates, translation, and/or rotation. Examples of such data include cases where each coordinate represents a fundamentally different feature (e.g.\ light, temperature, humidity). In Section \ref{sec:symmetrybreak}, we introduce two basic techniques, \emph{symmetry-breaking} and \emph{anchor points} to allow TDA to better study such data with heterogeneous features.

In this paper, we develop a novel {method} for topological machine learning for analyzing multivariate time series, with application to room occupancy detection. We use a dataset originating from the seminal paper by \citet{candanedo2016accurate}. In their research, data recorded from light, temperature, humidity and CO\textsubscript{2} sensors is provided. The main goal is to predict occupancy in an office room using these data. {We also include an additional experiment on Activity Recognition, using data from accelerometers.}

\begin{figure}[htbp]
\begin{center}
\begin{tikzpicture}[
    node distance = 5mm and 7mm,
      start chain = going right,
 disc/.style = {shape=cylinder, draw, shape aspect=0.3,
                shape border rotate=90,
                text width=20mm, align=center, font=\linespread{0.8}\selectfont},
  mdl/.style = {shape=ellipse, aspect=2.2, draw},
  alg/.style = {draw, align=center, font=\linespread{0.8}\selectfont}
                    ]
    \begin{scope}[every node/.append style={on chain, join=by -Stealth}]
\node (n1) [disc] {Multivariate\\ Time Series};
\node (n2) [alg] {Point\\ Clouds};
\node (n3) [alg] {Persistence\\ Diagrams};
\node (n4) [alg] {$k$-NN\\ (Wasserstein\\ distance)};
\node (n3) [mdl]  {Classification};
    \end{scope}
    \end{tikzpicture}
\end{center}
\caption{Outline of Topological Machine Learning for Multivariate Time Series.}
\label{fig:workflow}
\end{figure}
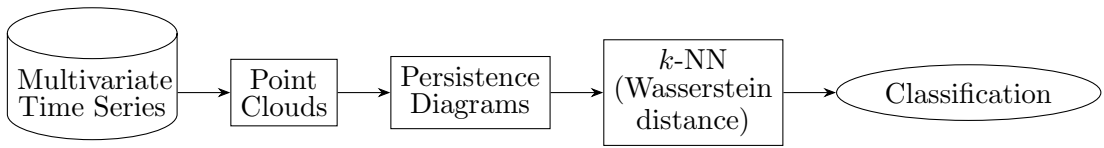

The outline of our {method} is summarized in Figure \ref{fig:workflow}. Firstly, we convert the multivariate time series to point cloud data via sliding windows \citep{gidea2018topological}. We also apply our techniques of symmetry-breaking and anchor points to the point clouds. Secondly, we generate persistence diagrams from the point cloud data. Lastly, we calculate the Wasserstein distance between the persistence diagrams and use the $k$-nearest neighbors algorithm ($k$-NN) for supervised machine learning (classification). {We will elaborate on each block in Figure} \ref{fig:workflow}, {both in theory and practice, in Sections} \ref{sec:methodology} {and} \ref{sec:experiments} {respectively.}

\subsection{Related work}
In the paper by \citet{gidea2018topological}, the authors introduce a new method based on topological data analysis to analyze financial time series and discover potential early signs of a financial crash. A key innovative factor in their paper is that their method can deal with multivariate time series (more than one time-dependent variable), which is different from the time-delay coordinate embedding for 1D time series \citep{takens1981detecting}. {Our paper generalizes their approach of converting a multivariate time series to a point cloud by introducing a new parameter representing stride. Other than that, our paper differs significantly in how we quantitatively study the data.} \citet{gidea2018topological} {makes use of the $L^p$-norm to study persistence landscapes, while our method uses the 1-Wasserstein distance of persistence diagrams to carry out supervised machine learning.}

In the paper by \citet{tran2019topological}, the authors study a delay-variant embedding method that constructs the topological features by considering the time delay as a variable parameter instead of considering it as a single fixed value. Their method studies multiple-time-scale patterns in a time series, which contains more information than just using a single time delay. {In their seminal work, a time series $x(t)$ is mapped to $m$-dimensional points using delay coordinates [$x(t),x(t-\tau),\dots,x(t-(m-1)\tau)]$ on the embedded space, where $\tau$ is a variable parameter denoting time delay and $m$ denotes the embedding dimension. Our proposed method differs from delay-variant embedding in two main ways. Firstly, our method works for multivariate time series, while the delay-variant embedding method focuses on univariate time series. Secondly, our proposed method did not consider time delay as a variable parameter.}

In the paper by \citet{merelli2016topological}, the authors study multivariate time series characterization using TDA, with applications to the epileptic brain. Their methodology is based on computing the Pearson correlation coefficients matrix for each window, followed by computing and plotting the weighted persistent entropy. {As a comparison, our method does not require the computation of the Pearson correlation coefficients matrix for each time window. Furthermore, our method includes a way of comparing quantitatively between two multivariate time series (via the Wasserstein distance of persistence diagrams). The paper by} \citet{merelli2016topological} {excels at comparing qualitatively two signals based on observing peaks in the plot of the respective weighted persistent entropies.}

\citet{umeda2017time} studied volatile time series using TDA. The methodology of the paper involves converting the time series into a quasi-attractor, and extracting topological information from the quasi-attractor in the form of a Betti sequence. In the learning step, a one-dimensional convolutional neural network (CNN) is used as a classifier. {For comparison, our method does not use CNN or other deep learning architectures. In addition, the paper assumes that an observed time series $\{x_1,\dots,x_t\}$ has a difference function $x_{k+1}=f(x_k,\dots,x_1)$ that specifies its behavior. For our method, we do not require or use this assumption.}

\citet{dirafzoon2016action} proposed a novel framework for activity recognition from 3D motion capture data using TDA. In the paper, point clouds are obtained from time series data using Takens' delay embedding \citep{takens1981detecting}. Subsequently, a feature vector for each time window is created using lengths of the most persistent off-diagonal features in the persistence diagram, with the maximum persistence interval length. As the final step, a nearest neighbor classifier with majority vote using Euclidean distance for the feature vectors is used, in order to classify each window. {We remark that our method is quite different in the way we use the topological information. Instead of a feature vector consisting of maximum persistence interval lengths, we use the Wasserstein distance of the persistence diagrams of each time window. In addition, the above paper uses preprocessing of the signals using a median filter for noise removal. For our method, no noise removal of the time series data was required.}

\citet{hofer2017deep} proposed and developed a technique that enables inputting topological signatures to deep neural networks and learn a task-optimal representation during training. An advantage of their method is that it learns the representation instead of mapping topological signatures to a pre-defined representation. {For comparison, our method does not use any deep learning technique or architecture. Also, the above paper mainly studies image (2D object shapes) and graph data, whereas we mainly focus on (multivariate) time series in this paper.}

In the paper by \citet{ravishanker2019topological}, the authors provide a comprehensive review of TDA for time series. In the work by \citet{seversky2016time}, the authors study the framework for the exploration of TDA techniques applied to time-series data. They consider and explore properties such as stability with respect to time series length, the approximation accuracy of sparse filtration methods, and the discriminating ability of persistence diagrams as a feature for learning. We note that both papers \citep{ravishanker2019topological,seversky2016time} utilize the time-delay coordinate embedding \citep{takens1981detecting}, also known as Takens' embedding \citep{mindlin1992topological}.

It is noted that due to the popularity and usefulness of time series in general, there are many other papers studying time series using TDA \citep{rucco2015topological,sanderson2017computational,gidea2018topological2,pita2017topological}.

We also remark that there are several established papers studying time series using other types of topology, in the broader sense of the word topology \citep{zhang2006complex,muldoon1993topology,tsonis2008topology,bonanno2003topology,djauhari2015optimality,mindlin1992topological}.

\subsection{Contribution}
Our paper combines the 4 key concepts: ``TDA'', ``machine learning ($k$-NN)'', ``multivariate'' and ``time series'', resulting in a novel {method} for topological machine learning on multivariate time series. Since TDA is a relatively new branch of data analysis, our paper also helps to validate and provide further evidence that topological methods work well in analyzing data. In addition, we demonstrate that TDA can be effectively combined with machine learning tools (e.g.\ $k$-NN algorithm) to study multivariate time series data. In addition, we also propose two basic methods, symmetry-breaking and anchor points, to study data that is sensitive to coordinates choice, translation and/or rotation. 

For applications, we demonstrate that our method can be effectively used to detect room occupancy. The detection of occupancy in buildings has been estimated to save energy in the order of 30\% to 42\% \citep{candanedo2016accurate,erickson2011observe,dong2009sensor}. Due to privacy concerns, it is also of interest to detect the presence of occupants without the use of a camera. Other applications for occupancy detection include security and analysis of building occupant behaviors \citep{candanedo2016accurate}.

{In addition, we apply our method on Activity Recognition (AR) data and obtained good results. We use data collected from accelerometers embedded on wearable sensors. AR is an emerging field of research with many potential applications such as monitoring the daily activity of the elderly (for ensuring their safety), fitness tracking and health informatics. There are also various other applications in the healthcare, human behavior modeling and human-machine interaction domains} \citep{palumbo2013multisensor}.

\section{Background}
We provide a brief overview of the relevant concepts in algebraic topology and persistent homology, and refer the reader to the appropriate sources for more details. A classical reference for algebraic topology is the text by \citet{Hatcher2002}, while the following papers provide an excellent introduction to persistent homology \citep{Ghrist2008,Zomorodian2005,edelsbrunner2008persistent}.

\subsection{Simplicial complexes}
\begin{figure}[!htbp]
\begin{center}
\begin{tikzpicture}
\node[circle,fill=black,inner sep=0pt,minimum size=3pt] (a) at (0,0){};
\node[below] at (0,-0.1) {$v_1$};
\node[below] at (0,-1) {vertex $\{v_1\}$};
\end{tikzpicture}
\begin{tikzpicture}
\node[circle,fill=black,inner sep=0pt,minimum size=3pt] (a) at (0,0){};
\node[circle,fill=black,inner sep=0pt,minimum size=3pt] (b) at (1.5,0){};
\node[below] at (0,-0.1) {$v_1$};
\node[below] at (1.5,-0.1) {$v_2$};
\node[below] at (0.75,-1) {edge $\{v_1,v_2\}$};
\draw (a) -- (b);
\end{tikzpicture}
\begin{tikzpicture}
\coordinate (a) at (0,0);
\coordinate (b) at (1.5,0);
\coordinate (c) at (0.75,1.3);
\filldraw[gray!50] (a.center) -- (b.center) -- (c.center) -- cycle;
\draw (a.center) -- (b.center) -- (c.center) -- cycle;
\node[circle,fill=black,inner sep=0pt,minimum size=3pt] at (a){};
\node[circle,fill=black,inner sep=0pt,minimum size=3pt] at (b){};
\node[circle,fill=black,inner sep=0pt,minimum size=3pt] at (c){};
\node[below] at (0,-0.1) {$v_1$};
\node[below] at (1.5,-0.1) {$v_2$};
\node[left] at (0.65,1.3) {$v_3$};
\node[below] at (0.75,-1) {triangle $\{v_1,v_2,v_3\}$};
\end{tikzpicture}
\begin{tikzpicture}
\coordinate (a) at (0,0);
\coordinate (b) at (1.5,0);
\coordinate (c) at (0.75,1.3);
\coordinate (d) at (1.8,0.6);
\begin{scope}[dashed,,opacity=0.5]
\draw (a) -- (d);
\end{scope}
\draw[fill=gray,opacity=0.5] (a)--(b)--(c);
\draw[fill=gray,opacity=0.5] (b)--(c)--(d);
\draw (a.center) -- (b.center) -- (c.center) -- cycle;
\draw (c) -- (d);
\draw (b) -- (d);
\node[circle,fill=black,inner sep=0pt,minimum size=3pt] at (a){};
\node[circle,fill=black,inner sep=0pt,minimum size=3pt] at (b){};
\node[circle,fill=black,inner sep=0pt,minimum size=3pt] at (c){};
\node[circle,fill=black,inner sep=0pt,minimum size=3pt] at (d){};
\node[below] at (0,-0.1) {$v_1$};
\node[below] at (1.5,-0.1) {$v_2$};
\node[left] at (0.65,1.3) {$v_3$};
\node[right] at (1.8,0.7) {$v_4$};
\node[below] at (0.5,-1) {tetrahedron $\{v_1,v_2,v_3,v_4\}$};
\end{tikzpicture}
\caption{A 0-simplex (vertex), 1-simplex (edge), 2-simplex (triangle) and 3-simplex (tetrahedron).}
\label{fig:simplices}
\end{center}
\end{figure}
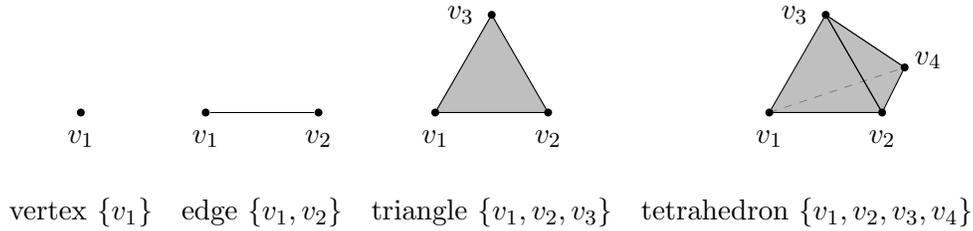

A \emph{simplicial complex} $K$ is a family of subsets of a set $S$ such that for every $\tau\subseteq\sigma\in K$, we have $\tau\in K$. The sets $\sigma\in K$ are called the \emph{faces} (or \emph{simplices}) of the simplicial complex $K$. We call the singleton sets $\{v\}$ the \emph{vertices} of $K$. The dimension of a simplex $\sigma\in K$ is defined to be $\dim(\sigma)=|\sigma|-1$, and we call a simplex of dimension $k$ a \emph{$k$-simplex}. Simplices of dimension 0, 1, 2, 3 can be viewed to represent a \emph{vertex}, \emph{edge}, \emph{triangle} and \emph{tetrahedron} respectively, as shown in Figure \ref{fig:simplices}.

A type of simplicial complex commonly used in TDA is the \emph{Vietoris-Rips complex} (or \emph{Rips complex} for short) which is defined as follows.

\begin{definition}
Let $\{x_i\}$ be a set of points in Euclidean space. The Rips complex $\mathcal{R}_\epsilon$ is the simplicial complex whose $k$-simplices are determined by each subset of $k+1$ points $\{x_j\}_{j=0}^k$ which are pairwise within distance $\epsilon$.
\end{definition}

We also introduce the concept of a \emph{filtration} of a simplicial complex $K$, which is a nested sequence of complexes $\emptyset=K^0\subseteq K^1\subseteq\dots\subseteq K^m=K$. We say that $K$ is a \emph{filtered complex}.

\subsection{Homology}
The $k$th \emph{chain group} $C_k$ of a simplicial complex $K$ is the free abelian group with basis the set of oriented $k$-simplices. The boundary operator $\partial_k: C_k\to C_{k-1}$ is defined on an oriented simplex $\sigma=[v_0,v_1,\dots,v_k]$ by
\begin{equation*}
\partial_k(\sigma)=\sum_{i=0}^k (-1)^i [v_0,\dots,\hat{v_i},\dots,v_k]
\end{equation*}
and extended linearly. The notation $\hat{v_i}$ indicates the deletion of the vertex $v_i$.

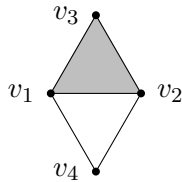
\begin{figure}[htbp]
\begin{center}
\begin{tikzpicture}[scale=0.8]
\coordinate (a) at (0,0);
\coordinate (b) at (1.5,0);
\coordinate (c) at (0.75,1.3);
\coordinate (d) at (0.75,-1.3);
\filldraw[gray!50] (a.center) -- (b.center) -- (c.center) -- cycle;
\draw (a.center) -- (b.center) -- (c.center) -- cycle;
\draw (a) --(d) --(b);
\node[circle,fill=black,inner sep=0pt,minimum size=3pt] at (a){};
\node[circle,fill=black,inner sep=0pt,minimum size=3pt] at (b){};
\node[circle,fill=black,inner sep=0pt,minimum size=3pt] at (c){};
\node[circle,fill=black,inner sep=0pt,minimum size=3pt] at (d){};
\node[left] at (-0.1,0) {$v_1$};
\node[right] at (1.6,0) {$v_2$};
\node[left] at (0.65,1.3) {$v_3$};
\node[left] at (0.65,-1.3) {$v_4$};
\end{tikzpicture}
\caption{For the above simplicial complex, we have $\beta_0=1$ (1 connected component), $\beta_1=1$ (1 circular hole which corresponds to the unshaded region) and $\beta_2=0$ (no ``voids'').}
\label{fig:betti}
\end{center}
\end{figure}

The \emph{cycle group} $Z_k$ and \emph{boundary group} $B_k$ are defined as $Z_k=\ker\partial_k$ and $B_k=\Ima\partial_{k+1}$ respectively. The $k$th \emph{homology group} is defined to be the quotient group $H_k=Z_k/B_k$. Informally, the rank of the $k$th homology group $\beta_k=\rank(H_k)$ (also called the $k$th Betti number) counts the number of $k$-dimensional holes in the simplicial complex $K$. For instance, $\beta_0$ counts the number of connected components (0-dim holes), $\beta_1$ counts the number of ``circular holes'' (1-dim holes), while $\beta_2$ counts the number of ``voids'' or ``cavities'' (2-dim holes). We show an example in Figure \ref{fig:betti}.

\subsection{Persistent homology}
Given a filtered complex $K$, the $i$th complex $K^i$ is naturally associated with the boundary operators $\partial_k^i$ and groups $C_k^i$, $Z_k^i$, $B_k^i$ and $H_k^i$. The \emph{$p$-persistent $k$th homology group} of $K^i$ is then defined as
\begin{equation*}
H_k^{i,p}=Z_k^i/(B_k^{i+p}\cap Z_k^i).
\end{equation*}

An equivalent definition of persistent homology groups is $H_k^{i,p}\cong\Ima \eta_k^{i,p}$, where $\eta_k^{i,p}: H_k^i\to H_k^{i+p}$ is the homomorphism that maps a homology class into the one that contains it \citep{ren2018weighted,Zomorodian2005}.

In brief, persistent homology studies a family of spaces parameterized by a distance $\epsilon$. The filtered complex $K$ is commonly obtained by the construction of Rips complexes over a range of distances $\epsilon$. Those topological features which persist over a parameter range can then be detected, revealing meaningful structures in the data.

\section{Methodology}
\label{sec:methodology}
In this section, we describe our methodology for studying multivariate time series using topological data analysis (TDA).

\subsection{Standardization of data}
\label{sec:standardize}
Following best practices in machine learning, we first standardize our data such that the values of each feature in the data have zero-mean and unit-variance. An advantage of standardizing is to prevent a feature with larger scale from completely dominating the other features.

\subsection{Converting multivariate time series to a point cloud}
\label{sec:pointcloud}

{The initial data type accepted for TDA is point cloud data. Hence, in order to study multivariate time series using topological methods, we would need to convert the multivariate time series to a point cloud. An example of such a conversion is the approach of }\citet{gidea2018topological}, {which we will describe in detail below}.

We will adopt the approach of \citet{gidea2018topological} to convert a multivariate time series to a point cloud data set, which is the required starting point for doing topological data analysis. We also generalize \citet{gidea2018topological} by introducing a new parameter $s$, representing stride.

We consider a multivariate time series consisting of $d$ 1-dimensional time series $\{x^k_n\}_n$, where $k=1,\dots,d$. Fix a sliding window of size $w$. For each time $t_n$, we define a point $x(t_n)=(x_n^1,\dots,x_n^d)\in\R^d$. Subsequently, for each time-window of size $w$, we obtain a point cloud 
\[X_n=(x(t_{1+s(n-1)}),x(t_{2+s(n-1)}),\dots,x(t_{w+s(n-1)}))\]
consisting of $w$ points in $\R^d$.

In brief, the length of the sliding window $w$ determines the size of the point cloud while the number $d$ of 1D time series determines the dimension of the point cloud \citep{gidea2018topological}. The stride $s$ determines how much the time-window slides for each consecutive point cloud. The value of $s$ corresponding to the original paper \citep{gidea2018topological} is a stride value of $s=1$.

In this paper, {for the first experiment on room occupancy}, we will choose a value of $w=s=10$, corresponding to non-overlapping sliding windows of length 10. {In the second experiment on activity recognition, we choose $w=s=5$, corresponding to non-overlapping windows of length 5. In principle, we can choose sliding windows of any length greater than 1. A sliding window of length 1 should be avoided as it leads to the point cloud having 1 point only (not counting the anchor point). A single point has trivial persistent homology and hence is not well suited for our topological method. By choosing different lengths of sliding windows in our two experiments, we demonstrate that our method can work for different lengths of sliding windows.}

\subsection{Symmetry-breaking and anchor points}
\label{sec:symmetrybreak}
In classical TDA, each coordinate plays the same role and has the same importance. For instance, in the case of $\R^3$, the $x$, $y$, $z$ coordinates are treated equally. Due to this symmetry property, topological methods excel in analyzing spatial data such as 3D point clouds \citep{singh2007topological,rosen2018inferring}.

However, this property may lead to TDA being unable to distinguish between certain point clouds. For example, persistent homology is unable to distinguish the two point clouds 
\begin{equation}
\label{eq:x1x2}
\begin{split}
X_1&=\{(0,0,0,0,0), (1,0,0,0,0)\},\\
X_2&=\{(0,0,0,0,0), (0,1,0,0,0)\},
\end{split}
\end{equation}
since in both point clouds the points are equidistant from each other (with distance 1). Alternatively, we can see that the two point clouds can be obtained from each other by rotation (and hence TDA is unable to distinguish them due to rotation-invariance). This can be a problem for certain data where each coordinate represents a fundamentally different type of feature (heterogeneous features). For instance, in our room occupancy data, the first coordinate represents temperature while the second represents humidity, so we \emph{do} actually want to distinguish between the two point clouds.

To this end, we introduce two basic techniques, symmetry-breaking and anchor points. Symmetry-breaking refers to adding a fixed constant vector to each point in the point cloud, while an anchor point refers to a fixed point that is introduced to the point cloud. We define them more precisely as follows.

\begin{definition}[Symmetry-breaking]
Let $X$ be a point cloud consisting of points in $\R^d$. Let $\mathbf{v}=(c_1,c_2,\dots,c_d)$ be a fixed vector in $\R^d$. We define the point cloud $X'$ obtained by \emph{symmetry-breaking} (of $X$) to be:
\begin{equation*}
X'=\{\mathbf{x}+\mathbf{v}\mid \mathbf{x}\in X\}.
\end{equation*}
\end{definition}

\begin{definition}[Anchor points]
Let $X$ be a point cloud consisting of points in $\R^d$. Let $A=\{\mathbf{a_1},\dots,\mathbf{a_n}\}$ be a set of points in $\R^d$, which we call \emph{anchor points}. 

We also define a new point cloud $Y=X\cup A$, called the point cloud augmented by the anchor points. 
\end{definition}

In this paper, we let $\mathbf{v}=(0,1,2,3,4)\in\R^5$ to be the fixed vector for symmetry-breaking. We take $A=\{(0,0,0,0,0)\}$, i.e., the only anchor point is the origin. With this choice, the point clouds in \eqref{eq:x1x2} become:

\begin{equation*}
\begin{split}
Y_1&=X_1'\cup A=\{(0,1,2,3,4),(1,1,2,3,4),(0,0,0,0,0)\}\\
Y_2&=X_2'\cup A=\{(0,1,2,3,4),(0,2,2,3,4),(0,0,0,0,0)\}.
\end{split}
\end{equation*}

We note that now the distance between $(1,1,2,3,4)$ and the origin $(0,0,0,0,0)$ is $\sqrt{31}$, while the distance between $(0,2,2,3,4)$ and the origin is $\sqrt{33}$. Hence, TDA is now able to distinguish between the point clouds $Y_1$ and $Y_2$ as desired.

\begin{figure}[!htbp]
\begin{center}
\begin{tikzpicture}
\draw[->,thick] (-2,0)--(2,0) node[right]{$x$};
\draw[->,thick] (0,-2)--(0,2) node[above]{$y$};
\node[circle,fill=black,inner sep=0pt,minimum size=5pt] (a) at (0,0){};
\node[left] at (-0.1,0.2) {anchor point};
\node[right] at (1.1,1) {point cloud};
\node[circle,fill=black,inner sep=0pt,minimum size=3pt] (a) at (1,1){};
\node[circle,fill=black,inner sep=0pt,minimum size=3pt] (a) at (0.9,0.9){};
\node[circle,fill=black,inner sep=0pt,minimum size=3pt] (a) at (0.75,0.9){};
\node[circle,fill=black,inner sep=0pt,minimum size=3pt] (a) at (0.9,0.75){};
\node[circle,fill=black,inner sep=0pt,minimum size=3pt] (a) at (0.75,0.7){};
\end{tikzpicture}
\begin{tikzpicture}
\draw[->,thick] (-2,0)--(2,0) node[right]{$x$};
\draw[->,thick] (0,-2)--(0,2) node[above]{$y$};
\node[circle,fill=black,inner sep=0pt,minimum size=5pt] (a) at (0,0){};
\node[left] at (-0.1,0.2) {anchor point};
\node[right] at (1.1,2) {point cloud (translated)};
\node[circle,fill=black,inner sep=0pt,minimum size=3pt] (a) at (1,2){};
\node[circle,fill=black,inner sep=0pt,minimum size=3pt] (a) at (0.9,1.9){};
\node[circle,fill=black,inner sep=0pt,minimum size=3pt] (a) at (0.75,1.9){};
\node[circle,fill=black,inner sep=0pt,minimum size=3pt] (a) at (0.9,1.75){};
\node[circle,fill=black,inner sep=0pt,minimum size=3pt] (a) at (0.75,1.7){};
\end{tikzpicture}
\caption{The inclusion of the anchor point at the origin can distinguish between the two point clouds which only differ by a translation. This is due to the differences in distance from the point clouds to the anchor point.}
\label{fig:translation}
\end{center}
\end{figure}
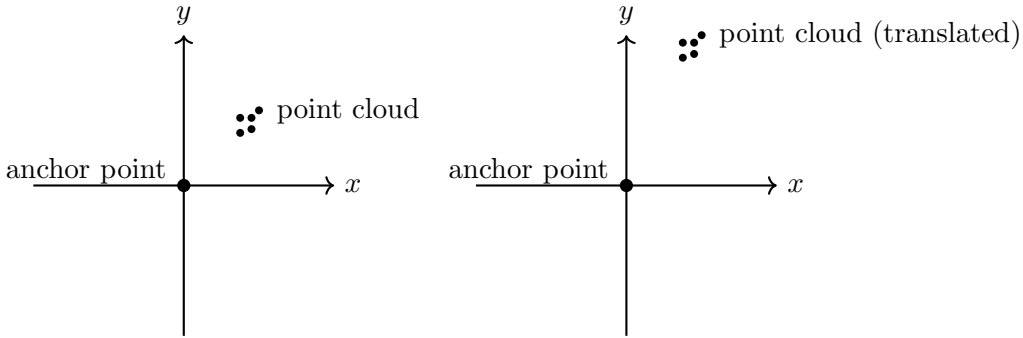

We also remark that the inclusion of anchor point(s) can further distinguish between point clouds that differ only by a translation. We illustrate this in Figure \ref{fig:translation}. Without the anchor point, TDA is generally unable to distinguish the two point clouds due to translation-invariance.

{We further summarize the proposed method of symmetry-breaking and anchor points in Algorithm} \ref{alg:algorithm}.

\begin{algorithm}[htbp]
    \caption{Symmetry-breaking and anchor points}\label{euclid}
    \hspace*{\algorithmicindent} \textbf{Input:} Point cloud $X$ consisting of points in $\R^d$.\\
    \hspace*{\algorithmicindent} \textbf{Output:} Augmented point cloud $Y=X'\cup A$.
    \begin{algorithmic}[1]
    \Procedure{SymmBreakAnchor}{$X$}
    \State $\mathbf{v} \gets (c_1,c_2,\dots,c_d)$
    \State $A \gets \{\mathbf{a_1},\dots,\mathbf{a_n}\}$
    \State $X' \gets \emptyset$
    \For{$\mathbf{x}$ in $X$}
    \State $X' \gets X'\cup\{\mathbf{x}+\mathbf{v}\}$
    \EndFor
    \State $Y\gets X'\cup A$
    \State \textbf{return} $Y$
    \EndProcedure
    \end{algorithmic}
    \label{alg:algorithm} 
\end{algorithm}

\subsection{From point cloud to persistence diagram}
A persistence diagram \citep{cohen2007stability} is a multiset of points in $\Delta:=\{(b,d)\in\R^2\mid b,d\geq 0, b\leq d\}$. Each point $(b,d)$ represents a generator of the homology group (of a chosen dimension), where $b$ denotes the birth of the generator and $d$ its death. In short, the persistence diagram can be viewed as a visual representation of the persistent homology of a point cloud. The persistence diagram is independent of choice of generators and thus is unique \citep{cohen2010lipschitz}. A key result is the stability of persistence diagrams with respect to Hausdorff distance, bottleneck distance \citep{cohen2007stability}, as well as Wasserstein distance \citep{cohen2010lipschitz}. Such stability results are desirable as it implies robustness against noise.

In this paper, we will focus on the Wasserstein distance with $p=1$, also known as the 1-Wasserstein distance or ``earth mover's distance''. The 1-Wasserstein distance is widely used in computer science to compare discrete distributions \citep{rubner2000earth,rabin2009statistical}. The Wasserstein distance is defined as follows \citep{cohen2010lipschitz,mileyko2011probability,berwald2018computing}.

\begin{definition}
The $p$-th Wasserstein distance between two persistence diagrams $D_1$, $D_2$ is defined as
\begin{equation*}
W_p(D_1,D_2)=\left(\inf_{\varphi:D_1\to D_2}\sum_{x\in D_1}\| x-\varphi(x)\|^p_\infty\right)^{1/p},
\end{equation*}
where the infimum is taken over all bijections $\varphi$ between $D_1$ and $D_2$.
\end{definition}

\subsection{The $k$-nearest neighbors algorithm}
To carry out classification (supervised machine learning), we utilize the $k$-nearest neighbors algorithm ($k$-NN) based on the Wasserstein distance. The $k$-NN algorithm is a relatively simple but yet effective machine learning algorithm that has been successfully applied across a wide range of domains \citep{batista2009k}.

For each point cloud $X$ (corresponding to a time window) in the test set, we will determine its $k$-nearest neighbors $\{Y_1,Y_2,\dots,Y_k\}$ in the training set, with respect to the Wasserstein distance. We then classify $X$ based on the majority class of the elements in the set $\{Y_1,Y_2,\dots,Y_k\}$.

\section{Experiments}
\label{sec:experiments}
The experiments were mainly implemented in Python, with the exception of computing the persistence diagrams and Wasserstein distances using the R package \texttt{TDA} \citep{fasy2014introduction}. The codes in the paper are made publicly available on GitHub: \url{https://github.com/wuchengyuan88/room-occupancy-topology}.

\subsection{Room occupancy}
\label{sec:roomoccupancy}

For this section, we study room occupancy data from the seminal paper by \citet{candanedo2016accurate}. In their setup, an office room with dimensions of 5.85m $\times$ 3.50m $\times$ 3.53m (W$\times$D$\times$H) was monitored for temperature, humidity, light and CO\textsubscript{2} levels. An Arduino microcontroller was used to acquire the data, and a digital camera was used to determine if the room was occupied or not. The data was recorded during the month of February (winter) in Mons, Belgium. The room was heated by hot water radiators \citep{candanedo2016accurate}.

Our experiment is regarding supervised machine learning (binary classification), where we train a model to predict if the room is non-occupied (class 0) or occupied (class 1) during a time window. We consider a room to be occupied if it is occupied during any period in the time window. That is, if a room is empty during some period in the time window, but occupied at other times in the same time window, we still classify it as occupied (class 1).

The 5 time-dependent variables in the data are the \texttt{Temperature}, \texttt{Humidity}, \texttt{Light}, \texttt{CO2}, \texttt{HumidityRatio} readings of each time period. Hence, the dimension of each point cloud is $d=5$. Measurements of each variable were taken after each time period of 1 minute. We divide the time periods into non-overlapping time windows of 10 minutes each (10 time periods). 

Following \citet{candanedo2016accurate}, we split the data into a training set and two test sets (Test Set 1 and Test Set 2), using a training--test ratio of 80:20. Following best practices in studying time series, we strictly respect the temporal order of the training and test sets. That is, our training and test sets come from distinct and non-overlapping time periods. Test Set 1 comes from time periods that are \emph{after} the training set, while Test Set 2 originates from time periods that are \emph{before} the training set. Benefits of having two such test sets include demonstrating that our method can predict \emph{future} as well as \emph{past} room occupancy (using the time-dependent variables). A further summary of the data sets can be found in Table \ref{table:datasetsummary}.

\begin{table}[!htp]
\footnotesize
\caption{Description of data sets.}
\begin{center}
\begin{tabular}{llll}
\hline
             &                        & \multicolumn{2}{l}{Data class distribution (\%)} \\\cline{3-4}
Data set     & Number of time windows & 0 (non-occupied)          & 1 (occupied)         \\ 
\hline
Training Set & 800                    & 77.5                      & 22.5                 \\
Test Set 1   & 200                    & 70.5                      & 29.5                 \\
Test Set 2   & 200                    & 57.0                        & 43.0                   \\
\hline
\end{tabular}
\end{center}
\label{table:datasetsummary}
\end{table}

\subsubsection{Standardization of data}
We standardize each of the 5 variables to zero-mean and unit-variance using StandardScaler from the Scikit-learn package. 

Technically, we are conducting two separate experiments, the first with Training Set and Test Set 1, and the second with Training Set and Test Set 2. Hence, we perform the standardization accordingly, first by standardizing the combined data set consisting of Training Set and Test Set 1, and then separately standardizing the combined data set consisting of Training Set and Test Set 2.

\subsubsection{Converting multivariate time series to a point cloud}
We follow the procedure outlined in Section \ref{sec:pointcloud}.

We set $w=s=10$ which corresponds to non-overlapping time windows of 10 minutes each. That is, each point cloud (not counting the anchor point) contains 10 points in $\R^5$.

\subsubsection{Symmetry-breaking and anchor points}
We use $\mathbf{v}=(0,1,2,3,4)$ as the fixed vector for symmetry-breaking. We use the origin $\mathbf{0}=(0,0,0,0,0)$ as the anchor point. That is, we augment each point cloud (corresponding to a time window) with the origin $\mathbf{0}$.

After symmetry-breaking, the mean and standard deviation (SD) of the 5 time-dependent variables for each experiment is as described in Table \ref{table:aftersymbreak}. 

\begin{table}[!htp]
\footnotesize
\caption{Mean and standard deviation (SD) of time-dependent variables.}
\begin{center}
\begin{tabular}{llllll}
\hline
     & Temperature & Humidity & Light & CO\textsubscript{2} & Humidity Ratio \\
     \hline
Mean & 0           & 1        & 2     & 3   & 4              \\
SD   & 1           & 1        & 1     & 1   & 1              \\
\hline
\end{tabular}
\end{center}
\label{table:aftersymbreak}
\end{table}

\subsubsection{From point cloud to persistence diagram}
\begin{figure}[!htbp]
\begin{center}
\includegraphics[scale=0.4]{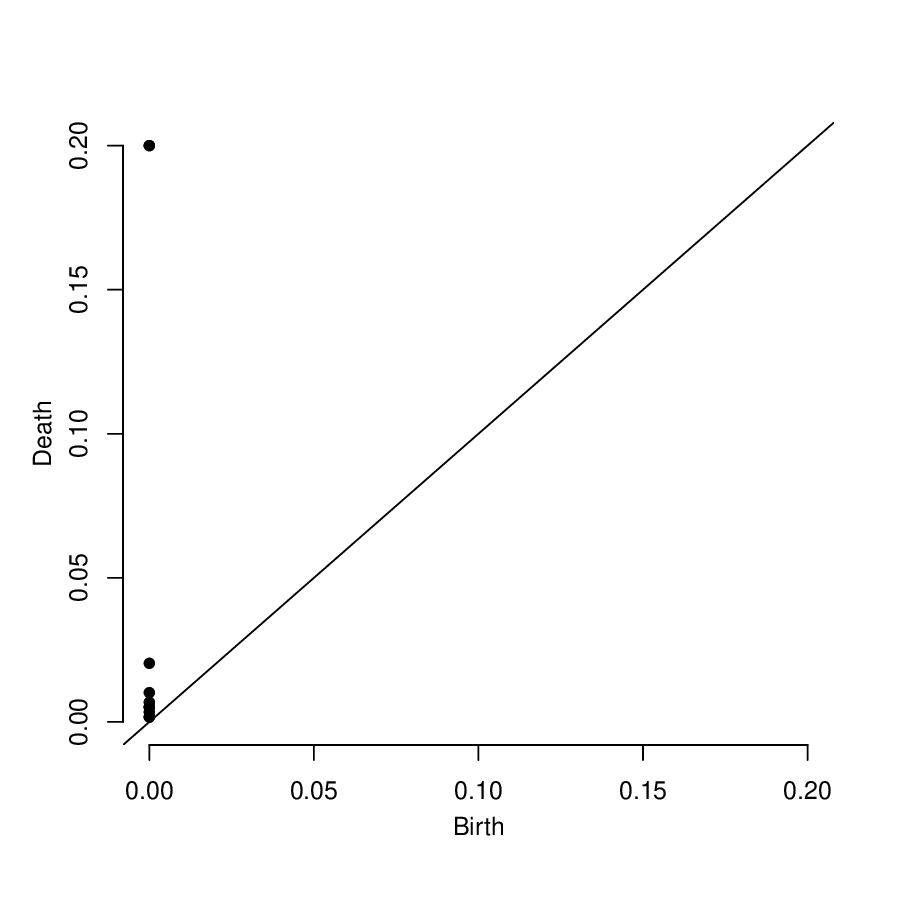}
\includegraphics[scale=0.4]{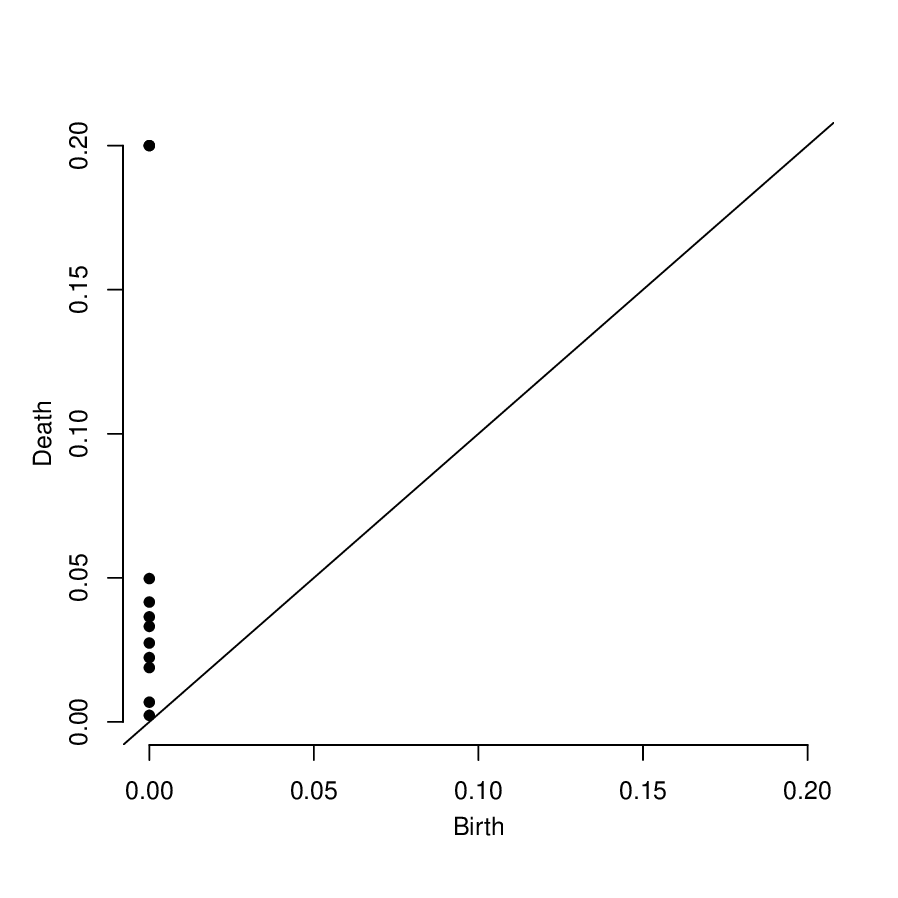}
\caption{The persistence diagram on the left belongs to a time window of class 0 (non-occupied), while that on the right belongs to a time window of class 1 (occupied). The points refer to homological features in dimension 0.}
\label{fig:twopd}
\end{center}
\end{figure}

For each point cloud, we construct the persistence diagram using the \texttt{ripsDiag} function in the R package \texttt{TDA}. The \texttt{ripsDiag} function uses a filtration of Rips complexes obtained from the point cloud to compute the persistence diagram. We show examples of two persistence diagrams from different classes in Figure \ref{fig:twopd}.

To calculate the 1-Wasserstein distance between two persistence diagrams, we use the \texttt{wasserstein} function, also from the R package \texttt{TDA}, with the default value of $p=1$. For this paper, we use the option \texttt{dimension=0} to specify that distances between persistence diagrams are computed using 0 dimensional features. This is because, for our data, we find that 1 dimensional (and higher) features rarely appear in the persistence diagrams, possibly due to the relatively small size of the point cloud.

\subsubsection{The $k$-nearest neighbors algorithm}
For this experiment, we choose $k=50$ for the $k$-NN algorithm. Since we have two test sets, we initially use Test Set 1 as a validation set to select a suitable value for the parameter $k$. Subsequently, we use the same $k$ for Test Set 2. We show the accuracy, sensitivity (true positive rate) and specificity (true negative rate) for various values of $k$ in Table \ref{table:differentk} (for Test Set 1). We select $k=50$ as the accuracy, sensitivity and specificity values are relatively good. We remark that we do not over-optimize for any single metric (e.g.\ accuracy) since other metrics are important as well in our study of room occupancy.

\begin{table}[!htp]
\footnotesize
\caption{Accuracy, sensitivity and specificity for different values of $k$ for Test Set 1.}
\begin{center}
\begin{tabular}{lllllllllll}
\hline
Value of $k$     & 10 & 20 & 30 & 40 & \textbf{50} & 60 & 70 & 80 & 90 & 100 \\
\hline
Accuracy (\%)    & 81 & 83 & 82 & 84 & \textbf{84} & 84 & 84 & 85 & 86 & 87  \\
Sensitivity (\%) & 88 & 88 & 80 & 81 & \textbf{80} & 76 & 69 & 68 & 64 & 63  \\
Specificity (\%) & 77 & 81 & 83 & 85 & \textbf{87} & 88 & 90 & 93 & 96 & 97 \\
\hline
\end{tabular}
\end{center}
\label{table:differentk}
\end{table}

To save computation time, we do not actually need to calculate the Wasserstein distance between all ${1000\choose 2} = 499 500$ pairs of persistence diagrams. We just need to calculate, for each of the 200 persistence diagrams in the test set, their respective Wasserstein distances from the 800 persistence diagrams in the training set. This amounts to $200\times 800=160000$ computations of Wasserstein distances, which is 68\% less computations than if we were to calculate distances between all pairs of diagrams.

The computation of Wasserstein distances is the most time-consuming part of the algorithm, but still taking a reasonably short time of 15 minutes (after the above reduction in computations) on a 2019 model of MacBook Pro with 2.4 GHz Intel Core i5 and 8 GB 2133 MHz LPDDR3.

\subsection{Activity Recognition}
{In addition, we perform a second experiment on Activity Recognition (AR) data. Our data is derived from the paper by} \citet{palumbo2013multisensor}, {which is available on the UCI Machine Learning Repository} \citep{Dua:2019}. {In their work, the authors collected data sampled by accelerometers embedded on wearable sensors as well as Received Signal Strength (RSS) of beacon packets exchanged between the sensors.}

{For this experiment, our goal is to predict the activity carried out by the user based on the given data. We focus on two activities: cycling (class 1) and standing up (class 0). We remark that the activity of standing up is an active one (not merely standing still), whereby there is a vertical momentum pattern of the standing up activity. The $y$-axis component of the accelerometer placed on the chest typically shows values indicative of the vertical direction of movement} \citep{palumbo2013multisensor}.

{The 6 time-dependent variables in the data are \texttt{avg\_rss12}, \texttt{var\_rss12}, \texttt{avg\_rss13}, \texttt{var\_rss13}, \texttt{avg\_rss23}, \texttt{var\_rss23}, where \texttt{avg} and \texttt{var} denote the mean and variance values of the RSS signals, respectively. Measurements were taken after each time period of 250 milliseconds. 
We divide the data into non-overlapping time windows of 5 time periods (1.25 seconds), where each time window represents a particular activity.}

{Subsequently, we split the data into training, validation and test sets using a ratio of 60:20:20. A summary of the data sets can be found in Table} \ref{table:datasetsummary2}.

\begin{table}[!htp]
\footnotesize
\caption{Description of data sets for activity recognition experiment.}
\begin{center}
\begin{tabular}{llll}
\hline
             &                        & \multicolumn{2}{l}{Data class distribution (\%)} \\\cline{3-4}
Data set     & Number of time windows & 0 (standing up)          & 1 (cycling)         \\ 
\hline
Training Set & 1728                    & 51.7                      & 48.3                 \\
Validation Set   & 576                    & 48.6                      & 51.4                 \\
Test Set    & 576                    & 46.4                        & 53.6                   \\
\hline
\end{tabular}
\end{center}
\label{table:datasetsummary2}
\end{table}

{The experiment is conducted in a similar manner as the first experiment in Section} \ref{sec:roomoccupancy}, {with two differences. Firstly, we choose $w=s=5$ corresponding to non-overlapping time windows of 5 time periods. The main reason for the above choice is to show that our method can work for time windows of different sizes. Secondly, due to there being 6 time-dependent variables, we use $\mathbf{v}=(0,1,2,3,4,5)$ as the fixed vector for symmetry-breaking.}

{For the $k$-NN algorithm, we use the validation set to select the optimal value of $k$. We select $k=40$ corresponding to a high validation accuracy of $98.61\%$, sensitivity of $99.32\%$ and specificity of $97.86\%$.}

\section{Results and discussion}
\subsection{Room Occupancy}
We obtain good results (80\% and above) for the key metrics of accuracy, sensitivity (recall of positive class) and specificity (recall of negative class) for both test sets. We summarize our results (including additional metrics such as precision and $F_1$ score) in Table \ref{table:finalresult}. {We remark that all metrics in Table }\ref{table:finalresult} {are computed by the \texttt{scikit-learn} package in Python. We consider these metrics as they are among the most popular in machine learning. Other than accuracy, metrics like sensitivity and specificity help to measure how an algorithm performs on an unbalanced dataset. We also state the confusion matrix $C$ where $C_{i,j}$ is the number of observations known to be in group $i$ and predicted to be in group $j$. The confusion matrices for Test Set 1 and Test Set 2 are}
$\begin{bmatrix}
122 & 19\\
12 &47
\end{bmatrix}$ {and}  
$\begin{bmatrix}
109 & 5\\
14 &72
\end{bmatrix}$
{respectively.}

{We also show a sample manual calculation, using information from the confusion matrix for Test Set 2, in Equation} \ref{eq:sample1}.

\begin{equation}
\label{eq:sample1}
\begin{split}
\text{Accuracy}&=\frac{\text{TP}+\text{TN}}{\text{TP}+\text{TN}+\text{FP}+\text{FN}}=\frac{72+109}{72+109+5+14}=0.91,\\
\text{Precision (class 1)}&=\frac{\text{TP}}{\text{TP}+\text{FP}}=\frac{72}{72+5}=0.94.
\end{split}
\end{equation}

\begin{table}[!htp]
\footnotesize
\caption{Results for Test Set 1 and Test Set 2 (using $k=50$).}
\begin{center}
\begin{tabular}{llllllll}
\hline
           &          &             &             & \multicolumn{2}{c}{Precision} & \multicolumn{2}{c}{$F_1$ score} \\\cline{5-6}\cline{7-8}
           & Accuracy & Sensitivity & Specificity & (class 0)     & (class 1)     & (class 0)      & (class 1)      \\
           \hline
Test Set 1 & 0.84     & 0.80        & 0.87        & 0.91          & 0.71          & 0.89           & 0.75           \\
Test Set 2 & 0.91     & 0.84        & 0.96        & 0.89          & 0.94          & 0.92           & 0.88          \\
\hline
\end{tabular}
\end{center}
\label{table:finalresult}
\end{table}

We remark that our results are not directly comparable with the seminal work by \citet{candanedo2016accurate}, even though we are using data derived from the same data set. This is because \citet{candanedo2016accurate} deals with real time occupancy detection (predicting room occupancy at each time period of 1 minute), while our work focuses on predicting room occupancy during a time window (of 10 time periods totalling 10 minutes).

There are also some additional difficulties in predicting room occupancy during a time window, especially with regards to correctly predicting non-occupancy (class 0). For instance, successfully predicting non-occupancy in a time window of length 10 is equivalent to 10 consecutive successful predictions of non-occupancy (for 10 time periods). Hence, even for a highly accurate model for real time prediction (e.g.\ 95\% accuracy), the chances of correctly predicting non-occupancy for a time window of length 10 drops to $0.95^{10}=59.9\%$ (assuming independence).

In view of the above discussions, our results show that our topological method is able to effectively and accurately predict room occupancy (and also non-occupancy) for time windows. Advantages of the time window approach include reducing the amount of data (many time periods are combined into a single time window) and hence computational time. In practice, energy saving measures also have good potential to work well with time windows. For example, it is practical to switch off the room air conditioner after the room has been empty for a time window, rather than immediately upon the room being vacant, since the occupants may only be temporarily leaving the room to return a short while later.

\subsection{Activity Recognition}
{We obtain very good results (above $99\%$) for all key metrics. We summarize our test set results in Table} \ref{table:finalresult2}. {The confusion matrix $C$ is }
$\begin{bmatrix}266 &1\\
0 &309
\end{bmatrix}$. {We remark that all metrics in Table }\ref{table:finalresult2} {are computed by the \texttt{scikit-learn} package in Python. We also show a sample manual calculation, using information from the confusion matrix, in Equation} \ref{eq:sample2}.

\begin{equation}
\label{eq:sample2}
\begin{split}
\text{Accuracy}&=\frac{\text{TP}+\text{TN}}{\text{TP}+\text{TN}+\text{FP}+\text{FN}}=\frac{309+266}{309+266+1+0}=0.9983,\\
\text{Specificity}&=\frac{\text{TN}}{\text{TN}+\text{FP}}=\frac{266}{266+1}=0.9963.
\end{split}
\end{equation}

\begin{table}[!htp]
\footnotesize
\caption{Results for Test Set.}
\begin{center}
\begin{tabular}{llllllll}
\hline
           &          &             &             & \multicolumn{2}{c}{Precision} & \multicolumn{2}{c}{$F_1$ score} \\\cline{5-6}\cline{7-8}
           & Accuracy & Sensitivity & Specificity & (class 0)     & (class 1)     & (class 0)      & (class 1)      \\
           \hline
  & 0.9983     & 1.0000        & 0.9963        & 1.0000          & 0.9968          & 0.9981           & 0.9984           \\
\hline
\end{tabular}
\end{center}
\label{table:finalresult2}
\end{table}

{Our good results are in line with those obtained in other previous works on activity recognition} \citep{palumbo2013multisensor,gallicchio2011user,lara2012centinela}. {Hence, it is a clear indication of the effectiveness of our proposed method for the activity recognition tasks of the type considered (cycling and standing up).}

{In addition, we remark that the time window chosen for this experiment is very short (1.25 seconds). Hence, this shows that our method can recognize and classify the activity types effectively even when given a very short time series consisting of data from the accelerometer sensors.}

\section{Conclusions}
This work provides a new {method} based on topological data analysis (TDA) to study multivariate time series. We use techniques in persistent homology (persistence diagram and Wasserstein distance) in combination with the $k$-nearest neighbors algorithm ($k$-NN) to perform supervised machine learning of time windows. In this paper, we also introduce methods (symmetry-breaking and anchor points) to allow TDA to better analyze data with heterogeneous features that are sensitive to translation, rotation, or choice of coordinates.

For applications, we first focus on room occupancy detection. Room occupancy detection is important in multiple ways, including energy saving, security and occupant behavior analysis. It is also important, for privacy reasons, to use non-intrusive types of data (such as light, humidity, temperature) instead of cameras (which may contain facial images of individuals) to detect room occupancy. Experimental results demonstrate the effectiveness of predicting room occupancy during a time window using topological methods.

{In the second application, we focus on activity recognition which is an emerging field of research. It has many potential applications such as infering the daily activity of the elderly for ensuring their safety, as well as applications in the healthcare, human behavior modeling and human-machine interaction domains} \citep{palumbo2013multisensor}. {The experimental results also demonstrate that topological methods are effective in recognizing activities based on accelerometer data from wearable devices.}

In a subsequent work, we also apply a variant of this topological method to analyze mixed numeric and categorical data \citep{wu2020topological}.

\section*{Disclosure statement}
No potential conflict of interest was reported by the authors.

\section*{Acknowledgements}
{The authors wish to thank the referees most warmly for numerous suggestions that have improved the exposition of this paper.}

\bibliographystyle{apacite}
\bibliography{bibtex11}
\end{document}